# WHAT DID FISHER MEAN BY AN ESTIMATE?

By Esa Uusipaikka[*]

*University of Turku*

Fisher's *Method of Maximum Likelihood* is shown to be a procedure for the construction of likelihood intervals or regions, instead of a procedure of point estimation. Based on Fisher's articles and books it is justified that by estimation Fisher meant the construction of likelihood intervals or regions from appropriate likelihood function and that an estimate is a statistic, that is, a function from a sample space to a parameter space such that the likelihood function obtained from the sampling distribution of the statistic at the observed value of the statistic is used to construct likelihood intervals or regions. Thus *Problem of Estimation* is how to choose the 'best' estimate. Fisher's solution for the problem of estimation is *Maximum Likelihood Estimate* (MLE). Fisher's *Theory of Statistical Estimation* is a chain of ideas used to justify MLE as the solution of the problem of estimation.

The construction of confidence intervals by the delta method from the asymptotic normal distribution of MLE is based on Fisher's ideas, but is against his 'logic of statistical inference'. Instead the construction of confidence intervals from the profile likelihood function of a given interest function of the parameter vector is considered as a solution more in line with Fisher's 'ideology'. A new method of calculation of profile likelihood-based confidence intervals for general smooth interest functions in general statistical models is considered.

**1. Introduction.** *'Collected Papers of R.A. Fisher'* (Bennet 1971) contains 294 articles. In eight of these (Fisher 1922, 1925b, 1932, 1934, 1935b, 1936, 1938, 1951) Fisher considers and explicitly mentions *"problem of estimation"*. In fact the title of the second of those articles is *'Theory of statistical estimation'*. Three of his books (Fisher 1925a, 1935a, 1956) all have a final chapter in which Fisher considers statistical estimation. In his last book *'Statistical Methods and Scientific Inference'* on page 143 Fisher, however, writes that

> A distinction without a difference has been introduced by certain writers who distinguish "Point estimation", meaning some process of arriving at an estimate without regard to its precision, from "Interval estimation" in which precision of the estimate is to some extent taken into account. "Point estimation" in this sense has never been practised either by myself, or by my

[*]Thanks to Professor Anthony Edwards for his valuable comments.

*Keywords and phrases:* confidence interval, delta method, estimate, likelihood interval, maximum likelihood of estimate, problem of estimation, profile likelihood, R. A. Fisher, theory of statistical estimation

predecessor Karl Pearson, who did consider the problem of estimation in some of its aspects, or by his predecessor Gauss of nearly one hundred years earlier, who laid the foundations of the subject.

In his famous R. A. Fisher Memorial Lecture (Savage 1976) in the year 1970 L. J. Savage reacts to this by saying that

> By "estimation," Fisher normally means what is ordinarily called point estimation. ... The term "point estimation" made Fisher nervous, because he associated it with estimation without regard for accuracy, which he regarded ridiculous and seemed to believe that some people advocated; . . .

and 26 years later B. Efron says in his R. A. Fisher Memorial Lecture (Efron 1998) that

> Fisher's great accomplishment was to provide an optimality standard for statistical estimation – a yardstick of the best it's possible to do in any given estimation problem. Moreover, he provided a practical method, maximum likelihood, that quite reliably produces estimators coming close to the ideal optimum even in small samples.

Even though Efron does not explicitly mention 'point estimation' later in his talk he speaks about maximum likelihood estimate and its estimated standard error, that is, about the well-known short hand expression

$$(1.1) \qquad\qquad\qquad \hat{\theta} \pm \text{se}_{\hat{\theta}}$$

(actually Efron considered approximate confidence intervals based on (1.1)). Savage's reaction may be interpreted as saying that hardly anybody uses point estimates without giving their estimated standard errors.

Fisher's comment on point estimation in Chapter *'The principles of estimation'* of his last book is puzzling because he himself uses in his articles and books heavily the notation (1.1) starting from his first article (Fisher 1912) on. The crucial word in that comment is 'precision'. On page 158 of the same chapter Fisher writes

> The study of the sampling errors, that is, of the precision, of statistical estimates, . . .

This quotation indicates that precision of an estimate meant for Fisher more than plain standard error.

In Section 2 Fisher's theory of statistical estimation is discussed. First concepts and notation used in this article will be introduced. Then a possible answer to the question in the title of this article, that is, to the question about the meaning of an estimate is given. This answer also shows that Fisher interpreted (1.1) differently and that difference explains his comment on point estimation and the difference of views between Fisher and all those who practice point estimation. Section 3 presents the prevailing view of maximum likelihood estimation (MLE) and shows that it is based on Fisher's ideas but is in conflict with his ideology of statistical inference. In Section 4 current state of Fisher's theory of estimation is reviewed. In the next section the special case of real-valued interest function is considered in more detail.

In writing this article two important decisions have been made. Firstly, because of the nature of the message of the article many quotations from Fisher's articles and books have been included. In connection of these additional citations to places which contain similar material are given. Secondly, technical material, for example, proofs of various propositions have not been included. These can be found in Fisher's publications and in the large literature discussing Fisher personally, his influence in statistics, and publications based on Fisher's scientific ideas. Historical material on Fisher includes Stigler (1973, 2005), Edwards (1974, 1997a, 1997b), Box (1978), Zabell (1989, 1992), and Aldrich (1997). Barnard (1963), Savage (1976), Rao (1992), and Efron (1998) are reviews about Fisher's influence on statistics. Early articles promoting likelihood-based inference include Bartlett (1936), Barnard (1948, 1951), Barnard et. al. (1962), Box and Cox (1964), and Kalbfleisch and Sprott (1970). Hacking (1965), Edwards (1972), Cox and Hinkley (1974), Barndorff-Nielsen (1988), Barndorff-Nielsen and Cox (1989, 1994), Lindsey (1996), Pace and Salvan (1997), Severini (2000), and Pawitan (2001) are books discussing the likelihood approach to statistical inference and contain further links to literature on Fisher.

## 2. Fisher's theory of statistical estimation.

### 2.1. *Statistical evidence.*

*Response and its statistical model.* Statistical inference is based on *statistical evidence*, which has at least two components. Two necessary components consist of the *observed response* $y_{\text{obs}}$ and its *statistical model* $\mathcal{M}$.

Because by definition the response contains random variation the actual observed value of the response is one among a set of plausible values for the response. The set of all possible values for the response is called *sample space* and denoted by $\mathcal{Y}$. The generic point $y$ in the sample space is called *response vector*.

The random variation contained in the response is described by either a point probability function or a density function defining a probability distribution in the sample space $\mathcal{Y}$. Because the purpose of statistical inference is to give information on the unknown features of the phenomenon under study these unknown features imply that the probability distribution of the response is not completely known. Thus it is assumed that there is a set of probability distributions capturing the important features of the phenomenon, among which features of primary interest are crucial, but also so-called secondary ones are important. It is always possible to index the set of possible probability distributions and in statistical inference this index is called *parameter*. The set of all possible values for the parameter is called the *parameter space* and denoted by $\Omega$. The generic point $\omega$ in the parameter space is called the *parameter vector*.

In this article only so-called *parametric statistical models* are considered. In parametric statistical models the value of the parameter is determined by the values of finitely many real numbers, which form the finite dimensional parameter vector. Thus the parameter space is a subset of some finite dimensional space $\mathcal{R}^d$ where $d$ is a positive integer denoting the dimension of the parameter space.

There is a probability distribution defined for every value of the parameter vector $\omega$. Thus a statistical *model function* is defined by a *model function* $p(y; \omega)$, which is a function of the observation vector $y$ and the parameter vector $\omega$ such that for every fixed value of parameter vector the function defines a distribution in the sample space $\mathcal{Y}$. The statistical model consists of the sample space $\mathcal{Y}$, the parameter space $\Omega$, and the model function $p(y; \omega)$. It is denoted by

$$(2.1) \qquad \mathcal{M} = \{p(y; \omega) : y \in \mathcal{Y}, \omega \in \Omega\}$$

or by the following mathematically more correct notation

$$(2.2) \qquad \mathcal{M} = \{p(\cdot; \omega) : \omega \in \Omega\},$$

where $p(\cdot; \omega)$ denotes the point probability or density function on the sample space $\mathcal{Y}$ defined by the parameter vector $\omega$.

*Statistical inference.*  Statistical inference concerns some characteristic or characteristics of the phenomenon from which the observations have arisen. The characteristics of the phenomenon under consideration are some functions of the parameter and are called the *parameter functions of interest*. If a subset of the components of the parameter vector form the interest functions, then the rest of the components are said to be *nuisance parameters*.

The result of a statistical inference procedure is a *collection of statements* concerning the unknown values of the parameter functions of interest. The statements are based on statistical evidence. The important problem of the *theory of statistical inference* is to characterize the form of statements that a given statistical evidence supports.

The *uncertainties* of the statements are also an essential part of the inference. The statements and their uncertainties are both based on statistical evidence and together form the *evidential meaning* of the statistical evidence.

*Likelihood concepts.* For the sake of generality let us denote by $A$ the event describing what has been observed with respect to the response. Usually $A$ consists just of one point of the sample space, but in certain applications $A$ is a larger subset of the sample space. The function $\Pr(A; \omega)$ of the parameter $\omega$, that is, the probability of the observed event $A$ with respect to the distribution of the statistical model is called the *likelihood function* of the model at the observed event $A$ and is denoted by $L_{\mathcal{M}}(\omega; A)$. If the statistical model consists of discrete distributions and the observed event has the form $A = \{\tilde{y}\}$, the likelihood function is simply the point probability function $p(\tilde{y}; \omega)$ as a function of $\omega$. If, however, the statistical model consists of absolutely continuous distributions, the measurement accuracy of the response has to be taken into account and so the event has the form $A = \{y : \tilde{y} - \delta < y < \tilde{y} + \delta\}$. Assuming for simplicity that the response is one-dimensional the likelihood function is

$$
\begin{aligned}
L_{\mathcal{M}}(\omega; \tilde{y}) &= \Pr(A; \omega) \\
&= \Pr(\tilde{y} - \delta < y < \tilde{y} + \delta; \omega) \\
&\approx p(\tilde{y}; \omega) 2\delta,
\end{aligned}
$$

which depends on $\omega$ only through the density function. Suppose that the response is transformed by a smooth one-to-one function $h$ and denote $\tilde{z} = h(\tilde{y})$. Let $g$ be the inverse of $h$. Then with respect to the transformed response

$$
\begin{aligned}
L_{\mathcal{M}}(\omega; \tilde{z}) &= \Pr(\tilde{z} - \delta < z < \tilde{z} + \delta) \\
&= \Pr(g(\tilde{z} - \delta) < y < g(\tilde{z} + \delta)) \\
&\approx p(\tilde{y}; \omega) 2 \left| g'(h(\tilde{y})) \right| \delta,
\end{aligned}
$$

which differs from the previous expression, but again depends on $\omega$ only through the original density function. From this follows that in the absolutely continuous case the likelihood function is up to a constant approximately equal to density function considered as the function of the parameter. The definition of the likelihood function which applies exactly to the discrete case and approximately to the absolute continuous case is

(2.3) $$ L_{\mathcal{M}}(\omega; \tilde{y}) = c(\tilde{y}) p(\tilde{y}; \omega), $$

where $c$ is an arbitrary positive function. The above discussion generalizes to models with higher dimensional parameter spaces. There are, however, situations in which the above approximation may fail (Lindsey 1996, pp. 75-80). In these situations it is wise to use the exact likelihood function, which is based on the actual measurement accuracy.

Occasionally the event $A$ is more complicated. For example, it might be known only that the response is greater than some given value $\tilde{y}$. Then $A = \{y : y > \tilde{y}\}$ and

$$L_{\mathcal{M}}(\omega; \tilde{y}) = \Pr(y > \tilde{y}; \omega) = 1 - F(\tilde{y}; \omega).$$

For theoretical reasons, which will be discussed later, instead of the likelihood function it is better to use *logarithmic likelihood function* or *log-likelihood function*, that is, the function

$$(2.4) \qquad\qquad l_{\mathcal{M}}(\omega; \tilde{y}) = \ln(L_{\mathcal{M}}(\tilde{y}; \omega)).$$

Often important information about the behavior of likelihood and log-likelihood functions can be obtained from first and second order derivatives of the log-likelihood function with respect to the components of the parameter vector. The vector of the first derivatives is called *score function* and negative of the matrix of second derivatives is called *observed information function*. Practically and theoretically the point in parameter space that maximizes the value of likelihood is of special importance.

Let $y_{\mathrm{obs}}$ be the observed response and $\mathcal{M} = \{p(y; \omega) : y \in \mathcal{Y}, \ \omega \in \Omega\}$ its statistical model. Likelihood ratio $L(\omega_1; y_{\mathrm{obs}})/L(\omega_2; y_{\mathrm{obs}})$ measures, how much more or less the observed response $y_{\mathrm{obs}}$ supports the value $\omega_1$ of the parameter vector compared to the value $\omega_2$. Because the likelihood ratio does not change, if the likelihood function is multiplied by a number not depending on the parameter vector $\omega$, it is convenient for the sake of various comparisons to multiply it by an appropriate number. A number generally used is the reciprocal of the maximum value of the likelihood function. The version of the likelihood function obtained in this way is called *relative likelihood function* and it has the form

$$(2.5) \qquad\qquad R(\omega; y_{\mathrm{obs}}) = \frac{L(\omega; y_{\mathrm{obs}})}{L(\hat{\omega}; y_{\mathrm{obs}})},$$

where $\hat{\omega}$ is the value of the parameter vector maximizing the likelihood function. The relative likelihood function takes values between 0 and 1 and its maximum value is 1. *Logarithmic relative likelihood function* or *relative*

Table 1

*Leukaemia data.*

| Treatment | Time of remission |
|-----------|-------------------|
| Drug | 6+, 6, 6, 6, 7, 9+, 10+, 10, 11+, 13, 16, 17+, 19, 20+, 22 |
| | 23, 25+, 32+, 32+, 34+, 35+ |
| Control | 1, 1, 2, 2, 3, 4, 4, 5, 5, 8, 8, 8, 8, 11, 11, 12,  12, 15, 17, 22, 23 |

Fig 1. *Leukemia data: Times of remission in leukemia data. Symbol* ○ *denotes times of the control group,* □ *uncensored times and* ■ *censored times of the drug group.*

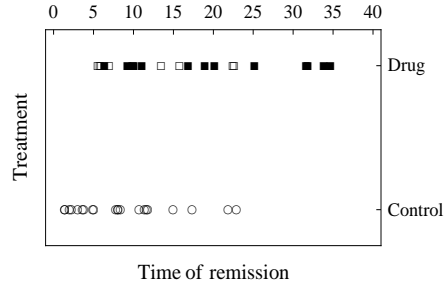

*log-likelihood function* is the logarithm of the relative likelihood function and thus has the form

$$(2.6) \qquad r(\omega; y_{\text{obs}}) = l(\omega; y_{\text{obs}}) - l(\hat{\omega}; y_{\text{obs}}).$$

The relative log-likelihood function has its values in the interval $(-\infty, 0)$ and its maximum value is 0.

*Example.* Data in Table 1 shows times of remission (i.e. freedom from symptoms in a precisely defined sense) of leukemia patients, some patients being treated with the drug 6-mercaptopurine (6-MP), the others serving as a control (Cox and Oakes 1984, p. 7). The columns of the data matrix contain values of treatment and time of remission. Censored times have been denoted by + sign. Figure 1 shows values of times of remission.

As an example consider the *'Drug'* group of the leukemia data. Assume that times of remission are a sample from some exponential distribution with unknown mean $\mu$. Statistical evidence consists of the response vector

$$y_{\text{obs}} = (1, 1, 2, 2, 3, 4, 4, 5, 5, 8, 8, 8, 8, 11, 11, 12, 12, 15, 17, 22, 23)$$

and statistical model

$$\mathcal{M} = \text{SamplingModel}[\text{ExponentialModel}[\mu], 21]$$

FIG 2. *Leukemia data: Relative likelihood function.*

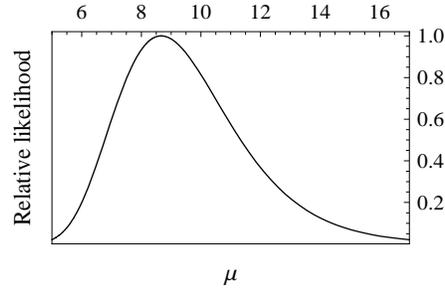

with sample space $\mathcal{Y} = \mathcal{R}^{21}$, parameter space $\Omega = (0, \infty)$, and model function

$$p(y; \mu) = \prod_{i=1}^{21} \frac{1}{\mu} e^{-\frac{y_i}{\mu}} = \frac{1}{\mu^{21}} e^{-\frac{21\bar{y}}{\mu}},$$

where $\bar{y}$ denotes the sample mean.

Likelihood function has the form

$$L(\mu; y_{\text{obs}}) = \frac{1}{\mu^{21}} e^{-\frac{182}{\mu}}$$

and log-likelihood function the form

$$l(\mu; y_{\text{obs}}) = -21\ln(\mu) - \frac{182}{\mu}.$$

Relative likelihood and relative log-likelihood functions are

$$R(\mu; y_{\text{obs}}) = \left(\frac{26}{3\mu}\right)^{21} e^{21 - \frac{182}{\mu}}$$

and

$$r(\mu; y_{\text{obs}}) = 21\ln\left(\frac{26}{3}\right) - 21\ln(\mu) + 21 - \frac{182}{\mu},$$

respectively. Figures 2 and 3 contains graphs of these relative likelihood functions.

2.2. *Problem of estimation.*

Fig 3. *Leukemia data: Relative log-likelihood function.*

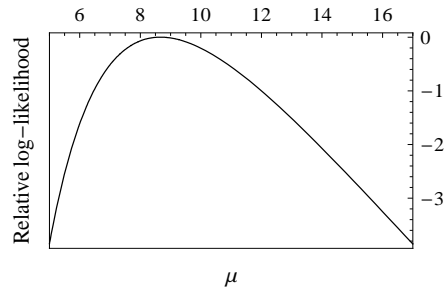

*Meaning of 'estimate'.* Fisher developed his theory of statistical estimation during 1912-56 by first inventing important concepts in his first four theoretical articles (Fisher 1912, 1915, 1920, 1921). In these articles, the 'problem of estimation' is not explicitly mentioned, but it will be argued that the problem and its solution appears in the last of these articles. In considering the development of Fisher's ideas and ideology it is important to separate concepts from the words Fisher used to name them. His terminology developed during these years and stabilized after ten years in Fisher (1922). The concepts, however, already appear during 1912-1921 in almost completely 'modern' form. Table 2 contains information on the development of concepts and their terminology for some important cases.

In Fisher (1956; see also 1921, p. 3; 1922, p. 313, 1925a, p. 8; 1925b, p. 701, 1935, p. 40) on page 49 the problem of estimation is defined.

> ... when the general hypothesis is found to be acceptable, and accepting it as true, we proceed to the next step of discussing the bearing of the observational record upon the problem of discriminating among the various possible values of the parameter, we are discussing the theory of estimation itself.

Already one should note that Fisher does not speak of picking one parameter value, but instead seems to be thinking of division of possible parameter values into sets.

In his first four articles (Fisher 1912, 1915, 1920, 1921) Fisher developed this solution by first introducing the concepts of likelihood and likelihood function in Fisher (1912). At the end of Fisher (1921, p. 25) he writes

> "Probable errors" attached to hypothetical quantities should not be interpreted as giving any information as to the probability that the quantity lies within any particular limits. When the sampling curves are normal and equivariant the "quartiles" obtained by adding and sub-tracting the probable error, express in reality the limits within which the likelihood exceeds 0.796,542, within twice, thrice, four times the probable error the values of the likelihood exceed 0.402,577, 0.129,098, and 0.026,267; within once, twice, and thrice the

TABLE 2
*Development of concepts and their terminology.*

| | | | | *Modern term* | | | |
|---|---|---|---|---|---|---|---|
| Year | *Parameter* | *Model function* | *Statistic* | *Sampling distribution* | *Likelihood* | *Likelihood function* | *MLE* |
| 1912 | arbitrary elements | function of known form | | frequency curve | inverse probability | inverse probability system | most probable set of values |
| 1915 | | frequency distribution of the population | statistical derivate | frequency distribution | | | most likely value |
| 1922 | parameter | hypothetical infinite population | statistic | distribution of statistics derived from samples | likelihood | | optimum value of parameter |

standard error, they exceed 0.606,051, 0.135,335 and 0.011,109.

Professor Anthony Edwards has checked the numbers and has noted that in case of the first standard error instead of 0.606 051 there should be 0.606 531. Other numbers are correct.

In Fisher (1956; see also 1922, p. 327; 1956, p. 53) on pages 69-72 a solution of the problem of estimation is given. On page 72 Fisher writes

> In the case under discussion a simple graph of the values of the Mathematical Likelihood expressed as a percentage of its maximum, against the possible values of the parameter $p$, shows clearly enough what values of the parameter have likelihoods comparable with the maximum, and outside of what limits the likelihood falls to levels at which the corresponding values of the parameter become implausible.

These quotations show that for Fisher the solution of the problem of estimation consisted of a collection of likelihood regions, that is, an inferential statement states that the unknown parameter (vector) belongs to a likelihood interval (region). Aldrich (1997) discusses the first quote and makes the same inference that Fisher was speaking about likelihood intervals.

This is, however, the solution only when the whole parameter vector is of interest and even in that case it lacks the all important assessment of uncertainty of the statement. In Fisher (1915, 1920, 1921) considered the case of real-valued interest functions, that is, correlation coefficient (1915, 1921) and standard deviation (1920). His solution in these articles was the construction of a statistic the sample distribution of which depended only on the interest function. The statistics chosen were naturally the sample correlation coefficients and the sample standard deviation.

At the beginning of his statistical writing Fisher wanted to express statistical inference, that is, statements concerning the unknown value of interest function in the from (1.1). There are at least three natural reasons for this, namely, first Fisher himself had learned statistics via Bayesian paradigm and so must have had strong inclination to use (1.1). Secondly this way of stating the inference should be easier for others to understand. Thirdly it gave a way to assess the uncertainty of the statement. Surely using (1.1) to express the inference was very problematic because people intended to think Fisher was just applying the Bayesian paradigm.

In addition to the problem of convincing others that he was not using the Bayesian paradigm, in which he did not succeed well, Fisher had a couple of other problems. First was the question of the interpretation of (1.1). Bayesian interpretation treats $\theta$ as random and $\hat{\theta}$ as fixed. Fisher instead treated $\theta$ as fixed and $\hat{\theta}$ as random. The second problem was that interval

$$(2.7) \qquad \hat{\theta} \pm k \, \mathrm{se}_{\hat{\theta}}$$

for some real number $k$ is a likelihood interval only if $\hat{\theta}$ has a normal distribution with constant standard deviation, that is, $\mathrm{se}_{\hat{\theta}}$ does not depend on $\theta$. Now, however, the invariance property of likelihood, which perhaps was the main motivation in Fisher (1912), came to rescue. Thus Fisher sought transformations $\psi = g(\theta)$ and $\hat{\psi} = g(\hat{\theta})$ such that $\hat{\psi}$ had (approximately) normal distribution with (approximately) constant standard deviation. Then he constructed a likelihood interval for $\psi$ using (2.7) and then transformed that to a likelihood interval for $\theta$, that is,

$$(2.8) \qquad g^{-1}(\hat{\psi} \pm k \, \mathrm{se}_{\hat{\psi}}).$$

At the same time, the problem of assessment of uncertainty was dealt by choosing the real number $k$ to be some quantile of standard normal distribution. It must be admitted, however, that Fisher in his earlier articles used the expression

$$(2.9) \qquad \hat{\psi} \pm \text{probable error of } \hat{\psi}$$

and later the expression

$$(2.10) \qquad \hat{\psi} \pm \text{standard error of } \hat{\psi}.$$

In his fourth article (Fisher 1921) on page 12 Fisher writes

It would of course be possible to render the statement of the possible error of $r$ less misleading by writing $\begin{smallmatrix} +0.1064 \\ -0.1361 \end{smallmatrix}$ instead of $\pm 0.1213$; that is by stating the actual quartile distances. Such a change, though certainly more accurate, and giving at any rate a danger signal as to the nature of the distribution, does not describe it effectively. Although two numbers are given, they contain less information than a single probable error when the distribution is normal; ...

That he interpreted (2.9) and (2.10) along with corresponding expressions for $\theta$ as collections of likelihood intervals he made clear at the end of Fisher (1921). That the expressions (2.9) and (2.10) implicitly contained an assessment of uncertainty of statements is evident.

Now the above explanation of what Fisher meant by the 'problem of estimation' and what was his solution to that problem explains his terminology, that is, the meaning he gave to the word 'estimate'. Thus for Fisher 'estimate', 'solution of estimation', and 'problem of estimation' meant the following.

**Estimate** consists of the observed values of some set of *statistics* which jointly define a function from sample space to parameter space.

**Solution of estimation** consists of determination of *likelihood regions* from the observed likelihood function based on the *sampling distribution* of those statistics.

**Problem of estimation** consists of deriving a method that produces best estimates.

So even though the observed value of an 'estimate' belonged to the parameter space that estimate for Fisher was not a point estimate. Instead this observed value of the estimate and its sampling distribution were carriers of information to be used in statistical inference. In the discussion of Savage (1976) Oscar Kempthorne expresses the same opinion

Savage alluded, appropriately, to obscurity on what Fisher meant by "estimation." My guess is that he meant the replacement of the data by a scalar statistic $T$ for the scalar parameter $\theta$ which contained as much as possible of the (Fisherian) information on $\theta$ in the data. But what one should do with an obtained $T$ was not clear, though Fisher was obviously not averse at times to regarding $T$ as an estimator of $\theta$. It is interesting, as Savage noted, that Fisher was the first to formulate the idea of exponential families in this connection. Here, also, the fascinating question of ancillaries arises, and on this Fisher was most obscure.

and similarly Ian Hacking in his book *'Logic of Statistical Inference'* (1965, p. 173).

But in Fisher's opinion, an estimate aims at being an accurate and extremely brief summary of the data bearing on the true value of some magnitude.

> Closeness on the true magnitude seems to be conceived as a kind of incidental feature of estimates.
>
> Thus from Fisher's point of view, he is quite correct when he writes 'if an unknown parameter $\theta$ is being estimated, any one-valued function of $\theta$ is necessarily being estimated by the same operation. The criteria used in the theory must, for this reason, be invariant for all such functional transformations of the parameters'. He is right, for by an estimate he means a summary of the data.
>
> It is often a trifle hard to see why what Fisher calls estimates should in any way resemble what are normally called estimates. And in fact his estimates based on small samples consist of a set of different numbers, only one of which resembles what we call an estimate, and the others of which are what he calls 'ancillary statistics', namely supplementary summaries of information.

Both Kempthorne and Hacking consider the case of one-dimensional parameter and in fact even though the above formulation of 'problem of estimation' includes also the case of a higher dimensional parameter it seems that Fisher never considered confidence regions in higher dimensional parameter spaces.

What Kempthorne considered 'obscure' is not that anymore if the word 'estimate' is interpreted the way presented above and Fisher's solution to the 'problem of estimation' which he gave at the end of Fisher (1921) is accepted. Even the 'question of ancillaries' has a natural explanation that will be considered in the next subsection.

*Solution of the problem of estimation.*  When comparing *mean error*

$$(2.11) \qquad \sqrt{\frac{\pi}{2}} \sum_{i=1}^{n} |y_i - \bar{y}|/n$$

with square root of *mean square error*

$$(2.12) \qquad \sqrt{\sum_{i=1}^{n} (y_i - \bar{y})^2/n}$$

as estimates (Fishers's meaning) of standard deviation in case of a sample from normal distribution with unknown mean and variance, Fisher hit, it may be said almost accidentally, on the concept of *sufficiency* (Fisher 1920, p. 768). The modern definition of sufficiency is the following.

**Sufficiency** A statistic $S$ is sufficient for the parameter vector $\omega$ if the conditional distribution of every other statistic $T$ with respect to $S$ is independent of the parameter vector $\omega$, that is, if the likelihood function obtained from the sampling distribution of $S$ is identical to the likelihood function obtained from the observed response and its statistical model.

For Fisher the term *sufficient statistic*, however, had a slightly different meaning. In addition to the term sufficient statistic Fisher used also the term *exhaustive statistic* and these had the following definitions.

**Sufficient statistic** is a statistic $S$ with dimension equal to that of the parameter vector $\omega$ such that the conditional distribution of every other statistic $T$ with respect to $S$ is independent of the parameter vector $\omega$, that is, a statistic $S$ such that the likelihood function obtained from the sampling distribution of $S$ is identical to the likelihood function obtained from the observed response and its statistical model.

**Exhaustive statistic** is a statistic $E$ which is sufficient in the modern meaning of the word such that it can be written in the form $E = (S, A)$ where $S$ has the same dimension as the parameter vector $\omega$ and the sampling distribution of $A$ is independent of the parameter vector $\omega$, that is, the likelihood function obtained from the conditional sampling distribution of $S$ given the observed value of $A$ is identical to the likelihood function obtained from observed response and its statistical model.

Thus for Fisher a $d$-dimensional statistic $S$ is sufficient if it is sufficient in the modern meaning of the word. Now, because Fisher's estimate meant a $d$-dimensional statistic taking values from the parameter space, he used in fact the term *sufficient estimate*.

In both cases Fisher's aim was an estimate (Fisher's meaning) such that the likelihood function obtained from the conditional sampling distribution of the estimate will exhaust all the information on the parameter vector contained in the statistical evidence. In case of Fisher's sufficiency there is no need to condition and so the estimate is sufficient because it exhausts all the information.

In Fisher (1920) it was shown that square root of mean square error is the sufficient estimate for standard deviation.

2.3. *Theory of statistical estimation.*

*Background.* Fisher (1922) is the first of Fisher's articles that considered the problem of estimation in a systematic way and contains first version of a chain of ideas that Fisher called *Theory of statistical estimation* and which was meant to give a justification for his solution of the problem of estimation. The previous four articles (1912, 1915, 1920, 1921) already contained, however, seeds of these ideas and especially four basic 'principles' that form Fisher's 'ideology' which he kept talking about for the rest of his life. The principles were the following.

**No prior distribution** Statistical inference on unknown parameter(s) should be done without assuming any prior distribution for the unknown parameters, unless the prior distribution was based on real physical knowledge of the situation.

**Invariance** Statistical inference should be invariant with respect one-to-one transformations of parameters.

**Efficiency** Statistical inference should be efficient so that all the information in the observations should be used in an 'optimal' way.

**Small-samples** Statistical inference should be exactly applicable in small-samples.

In early papers, these principles were more Fisher's reactions against certain common ways of thinking. The first principle was based on the writings of Boole (1854), Venn (1866), etc. and expresses Fisher's conviction that "the theory of inverse probability is founded upon an error, and must be wholly rejected (Fisher 1925a, p. 9)." Here the term *inverse probability* refers to Laplacean way of using 'uniform' prior distributions when there is no knowledge about the parameters. It must, however, be said that Fisher did not reject Bayesian inference altogether. On the contrary, he regarded it highly, but insisted that it was applicable only when prior distribution could be based on real physical knowledge and not on the lack of it.

Fisher's favorable thoughts about Bayes and Bayesian inference have an explanation in his ideas about the uncertainty of statistical inference. In Fisher (1956, p. 40) he writes the following.

> While, as Bayes perceived, the concept of Mathematical Probability affords a means, in some cases, of expressing inferences from observational data, involving a degree of uncertainty, and of expressing them rigorously, in that the nature and extent of the uncertainty is specified with exactitude, yet it is by no means axiomatic that the appropriate inferences, though in all cases involving uncertainty, should always be rigorously expressible in terms of this same concept.

In this quotation the term *Mathematical Probability* means a probability model, that is, a sample space of outcomes, a collection of events, and a probability measure defined on the collection of events giving the distribution. As the quotation indicates according to Fisher there are different situations which afford different kinds of measures of uncertainty. Mathematical Probability is on the top of these measures of uncertainty. This explains Fisher's positive attitude to Bayesian inference and his *theory of fiducial probability* which is an attempt to produce statistical inferences in the form of Mathematical Probability without assuming prior distribution for parameters. Other forms of uncertainty include *Mathematical Likelihood*, significance levels, and confidence levels. In these latter cases, the calculated numbers are measures of uncertainties of the conclusions made, but they cannot be expressed in the form of Mathematical Probability because there is no collection of events and probability measure on it.

The principle of invariance appeared already in Fisher (1912) and in Fisher (1956, p. 146) it is said that

> ...if an unknown parameter $\theta$ is being estimated, any one-valued function of $\theta$ is necessarily being estimated by the same operation. The criteria used in the theory must, for this reason, be invariant for all such functional transformations of the parameters.

This principle has been the most important from the beginning of Fisher's statistical writing. The word 'absolute' appearing in the title of Fisher (1912) refers to invariance (Aldrich 1997 and Stigler 2005). Fisher used the concept of invariance also to downplay the importance of concepts like unbiasedness, etc.

One of the reasons of writing Fisher (1912) was to present an alternative to the method of moments and to show that method of moments was not a method that generally produced efficient solutions. In Fisher (1922) he gave a dramatic example of the failure of method of moments by showing that the use of the sample mean calculated from the sample from Cauchy distribution with unknown location and known scale was equivalent to basing the inference on just one observation and discarding all the other observations in the sample.

From the beginning of his statistical writing Fisher wanted to develop solutions to small-sample situations. Before Fisher came to the scene, most methods of statistical inference were applicable when samples were large. The first paper on small-samples was Student (1908). Fisher admired the work of Gosset and the influence of this work was Fisher's emphasis on finite samples and his view that the duty of statistics and statisticians is to provide exact methods of statistical inference that scientist can use to analyze their data, even when that data consists of a small sample.

*Fisher's logic of inductive inference.* Fisher's *Theory of statistical estimation* is a chain of ideas intended to justify his solution of the problem of estimation. In Fisher (1935, p. 41) he thought it necessary to

> ...show how it is that a consideration of the problem of estimation, without postulating any special significance for the likelihood function, and of course without introducing any such postulate as that needed for inverse probability, does really demonstrate the adequacy of the concept of likelihood for inductive reasoning, in the particular logical situation for which it has been introduced ...

Stigler (2005) contains a thorough discussion about the developments that led to Fisher (1922) in which Fisher presented first version of theory of statistical estimation. Stigler's view about the motive behind the theory of statistical estimation is perfectly in line with views presented in this article.

The chain of ideas in Fisher's *Theory of statistical estimation* naturally divides in two parts. In Fisher (1935, p. 41) he describes this chain in the following way.

> This logical characteristic of our approach naturally requires that our edifice shall be built in two stories. In the first we are concerned with the theory of theory *large samples*, using this term, as is usual, to mean that nothing that we say shall be true, execpt in the limit when the size of the sample is indefinitely increased; a limit, obviously, never attained in practice. This part of the theory, to set off against the complete unreality of its subject-matter, exploits the advantage that in this unreal world all possible merits of an estimate may be judged exclusively from its variablity, or sampling variance. In the second story, where the *real* problem of finite samples is considered, the requirement that our estimates from these samples may be wanted as materials for a subsequent process of estimation is found to supply the unequivocal criteria required.

In this quotation, italics come from Fisher and is important because it emphasizes Fisher's often expressed opinion that ideas of the first 'story' are theoretical and were not intended to be used in finite samples.

In the first story there appears three concepts, namely, *consistency*, *efficiency*, and *expected information* or *Fisher information*. Consistency means that in the limit estimate becomes constant that is equal to the parameter function of interest. Efficiency means that the asymptotic variance of the estimate is as small as possible. Technically, as 'sample size' $n$ increases the limiting value of $nV$, where $V$ stands for the variance of our estimate, shall be as small as possible (Fisher 1935, p. 42). Fisher showed that for any consistent estimate

$$V \geq \frac{1}{i},$$

where

$$i = \mathrm{var}\left(\frac{\partial l(\omega)}{\partial \omega}\right)$$

is the Fisher information. Fisher showed that MLE is consistent and its limiting variance is equal to the reciprocal of Fisher information. Thus in the limit MLE is the best, but this result concerns what happens in an unreal world.

One may ask why Fisher included in his theory of statistical estimation the first story of large samples. An explanation might be that he wanted to show the connection between his theory and the standard interval (1.1). It is well-known that in most cases the asymptotic distribution of an estimate $\hat{\theta}$ is normal. In order that (1.1) would be a likelihood interval for $\theta$ the mean of the asymptotic normal distribution should be exactly or at least approximately equal to $\theta$ and variance of it should be exactly or approximately independent of $\theta$. The first requirement is true if the estimate is consistent. The other reason may be that the concept of Fisher information, which was introduced in the first story and had there a rather simple interpretation, was applicable also in the second story, that is, in small samples.

The first 'story' in Fisher's theory of statistical estimation was just a prelude to the discussion of the real problem in the second 'story'. In Fisher (1935, p. 46; see also 1925a, p. 15; 1925b, p. 712; 1956, p. 46, 158, 163) he says the following.

> We are now in a position to consider the real problem of finite samples. For any method of estimation has its own characteristic distribution of errors, not now necessarily normal, and therefore its own intrinsic accuracy.

Because in finite samples the sampling distributions of the various possible estimates can have widely different forms and can no more be compared with the help of standard deviation Fisher needed some new criterion for the comparison of 'efficiency' of different estimates. Early on he noticed that the concept of Fisher information was general and can be determined in all situations. So in the second 'story', that is, in the real problem of finite samples Fisher information replaced asymptotic variance as a tool used to rank estimates. In Fisher (1935, p. 46; see also 1925a, p. 314, 338; 1925b, p. 709, 712, 714; 1956, p. 153, 157) he wrote

> This quantity $i$, which is independent of our methods of estimation, evidently deserves careful consideration as an intrinsic property of the population sampled. In the particular case of error curves, or distributions of estimates of the same parameter, the amount of information of a single observation evidently provides a measure of the intrinsic accuracy with which it is possible to evaluate that parameter, and so provides a basis for comparing the accuracy of error curves which are not normal, but may be of quite different forms.

Fisher showed that

$$0 \leq i_{\tilde{\theta}}(\theta) \leq i(\theta),$$

where $\tilde{\theta}$ is some estimate of $\theta$. Also $i(\theta)$ is the Fisher information calculated from the original data and $i_{\hat{\theta}}(\theta)$ is that calculated from the sampling distribution of the estimate $\tilde{\theta}$. When $\hat{\theta}$ is a sufficient estimate its Fisher information equals that of original data. In this case solution of the problem of estimation consists of the likelihood function obtained from the sampling distribution of the maximum likelihood estimate evaluated at the observed value of MLE. But from sufficiency follows that this likelihood function is just equivalent to the likelihood function obtained from the original data. In Fisher (1935, p. 47) it is said that

> Having obtained a criterion for judging the merits of an estimate in the real case of finite samples, the important fact emerges that, though sometimes the best estimate we can make exhausts the information in the sample, and is equivalent for all future purposes to the original data, yet sometimes it fails to do so, but leaves a measurable amount of the information unutilized. How can we supplement our estimate so as to utilize this too? It is shown that some, or sometimes all the lost information may be recovered by calculating what

> I call ancillary statistics, which themselves tell us nothing about the value of
> the parameter, but instead, tell us how good an estimate we have made of it.

In (Fisher,1922) Fisher thought first that MLE is always sufficient estimate, but realized almost immediately this is not true. In this case

$$i_{\hat{\theta}}(\theta) < i(\theta),$$

and the likelihood function constructed from the sampling distribution MLE does not contain all the relevant information about the unknown parameter. So Fisher had a serious problem of how the lost information could be recovered. In (Fisher, 1934) he found two important special cases in which the lost information could be recovered by considering the conditional sampling distribution of MLE given the observed values of certain ancillary statistics, that is, statistics the marginal distribution of which did not depend on the unknown parameter. He showed that the expected value of the conditional Fisher information calculated from the conditional sampling distribution of MLE is equal to the Fisher information of the original data. Now also in this case the solution of the problem of estimation consists just of the original likelihood function because the likelihood function obtained from the conditional sampling distribution of MLE is equal to the original likelihood function.

The net result of Fisher's *Theory of statistical estimation* is that solution of the problem of estimation consists of likelihood intervals obtained from the likelihood function of data. After discussing the role ancillary statistic in (Fisher 1956, p. 161) Fisher gives the following summary.

> ...it is the Likelihood function that must supply all the material for estimation, and that the ancillary statistics obtained by differentiating this function are inadequate only because they do not specify the function fully.

2.4. *Example.* As an example, Fisher (1925, p. 705) contained the comparison of sample mean and median as estimates in the sense Fisher meant by estimates. Assume that the response vector is a sample of size $n = 11$ from some normal distribution with unknown mean $\mu$ and known standard deviation $\sigma = 1$. The log-likelihood function $l_1$ calculated from the sampling distribution of the sample mean $\bar{y}$ has the form

$$l_1(\mu) = -n(\bar{y} - \mu)^2/2$$

and the log-likelihood function $l_2$ calculated from the sampling distribution of the sample median $\tilde{y}$ has the form

$$l_2(\mu) = -(\tilde{y} - \mu)^2/2 + (n-1)\ln(\Phi(\tilde{y} - \mu))/2 + (n-1)\ln(1 - \Phi(\tilde{y} - \mu))/2.$$

Fig 4. *Log-likelihood functions of mean obtained from sampling distribution of sample mean (continuous curve) and sample median (dashed curve).*

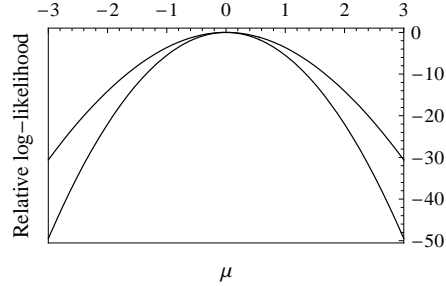

Because here the purpose is to compare the shape of these log-likelihood functions assume without restriction that $\bar{y} = \tilde{y} = 0$. Figure 4 shows the graphs of these log-likelihood functions and from these it is clearly seen that the sample mean contains more information on $\mu$ than the sample median. In fact, as is well-known, the sample mean is in this model sufficient statistic, both in Fisher's and modern meaning, and thus contains all the information on $\mu$ there is in the statistical evidence. The sample median instead loses some of this information and is an inefficient estimate. The asymptotic variances of the estimates are $\sigma^2/n$ and $\pi\sigma^2/2n$, respectively. The asymptotic efficiency of the sample median is thus $2/\pi = 63.33\%$.

## 3. Prevailing anti-fisherian view of MLE or Maximum Likelihood Estimate Method.

3.1. *Delta method and Wald confidence interval.* Let $g(\omega)$ be a given real valued interest function with $\psi$ as value. The most often used procedure to construct a confidence interval for $\psi$ is so-called *delta-method*. The delta-method interval has the form

$$(3.1) \qquad \psi \in \hat{\psi} \mp z^*_{\alpha/2} \sqrt{\frac{\partial g(\hat{\omega})}{\partial \omega}^{\mathrm{T}} J(\hat{\omega})^{-1} \frac{\partial g(\hat{\omega})}{\partial \omega}},$$

where $\hat{\psi} = g(\hat{\omega})$ is the maximum likelihood estimate of $\psi$, $J(\omega)$ the observed information matrix, and $z^*_{\alpha/2}$ the $(1 - \alpha/2)$-quantile of standard normal distribution. When interest function consists of a component $\omega_i$ of the parameter vector $\omega$, the above interval takes the form

$$(3.2) \qquad \omega_i \in \hat{\omega}_i \mp z^*_{\alpha/2} s_{\hat{\omega}_i},$$

where the *standard error of estimate*

$$(3.3) \qquad s_{\hat{\omega}_i} = \sqrt{(J(\hat{\omega})^{-1})_{ii}}$$

is the square root of the $i^{\text{th}}$ diagonal element of the inverse $J(\hat{\omega})^{-1}$. This latter confidence interval is better known by the name *Wald confidence interval*.

Wald confidence interval is derived from the asymptotic normal distribution of the maximum likelihood estimate $\hat{\omega}$, that is,

$$(3.4) \qquad \hat{\omega} \simeq \text{MultivariateNormal}[d, \omega, I(\omega)^{-1}],$$

where $I(\omega)$ is the Fisher expected information matrix which can be estimated by the observed information matrix $J(\hat{\omega})$. Thus first we have from the asymptotic normal distribution

$$\frac{\hat{\omega}_i - \omega_i}{\sqrt{(I(\omega)^{-1})_{ii}}} \simeq \text{NormalModel}[0, 1]$$

and then by inserting the estimate

$$\frac{\hat{\omega}_i - \omega_i}{\sqrt{(J(\hat{\omega})^{-1})_{ii}}} \simeq \text{NormalModel}[0, 1].$$

from which the Wald confidence interval follows.

The delta-method interval is then obtained from the Taylor-series expansion

$$\hat{\psi} = g(\hat{\omega}) \approx g(\omega) + \frac{\partial g(\omega)}{\partial \omega}^{\text{T}} (\hat{\omega} - \omega),$$

which gives

$$g(\hat{\omega}) \simeq \text{NormalModel}\left[\psi, \frac{\partial g(\omega)}{\partial \omega}^{\text{T}} I(\omega)^{-1} \frac{\partial g(\omega)}{\partial \omega}\right].$$

So

$$\frac{\hat{\psi} - \psi}{\sqrt{\frac{\partial g(\omega)}{\partial \omega}^{\text{T}} I(\omega)^{-1} \frac{\partial g(\omega)}{\partial \omega}}} \simeq \text{NormalModel}[0, 1],$$

from which by replacing the denominator by its estimator we get

$$\frac{\hat{\psi} - \psi}{\sqrt{\frac{\partial g(\hat{\omega})}{\partial \omega}^{\text{T}} J(\hat{\omega})^{-1} \frac{\partial g(\hat{\omega})}{\partial \omega}}} \simeq \text{NormalModel}[0, 1].$$

It can be shown that Wald and delta-method intervals are approximations of profile likelihood-based confidence intervals. These approximations have, however, some serious drawbacks, the most serious of which is that they are not invariant with respect to monotone transformations of the interest function. Secondly, these intervals often include impossible values of the interest function.

3.2. *MLE-based and likelihood-based confidence intervals.* The following quotations from books which have been written by men who promote likelihood approach show general view on Fisher's 'Method of Maximum Likelihood'. In most other theoretical books and articles on maximum likelihood estimation, the approach is similar and even more uncompromising in the sense that they do not contain reservations expressed by Anthony Edwards, James Lindsey, and Yudi Pawitan.

Edwards writes in *Likelihood* (1972, p. 98)

> He [Fisher] advocated what he later called the Method of Maximum Likelihood in his very first paper, as a means of point estimation.

Lindsey in *Parametric Statistical Inference* (1996, p. 81)

> The maximum likelihood estimate can be looked upon as a point estimate. However, like any point estimate, in most contexts, outside of those just mentioned, it is often of little use because many other models will be almost as likely. We need to look at the form of the whole likelihood function, . . .

and Pawitan in *In All Likelihood: Statistical Modelling and Inference Using Likelihood* (2001, p. 30)

> Fisher (1922) introduced likelihood in the context of estimation via the method of maximum likelihood, but in his later years he did not think of it as simply a device to produce parameter estimates.

Contrary to what Edwards, Lindsey, and Pawitan say by the term 'Method of Maximum Likelihood' Fisher meant that maximum likelihood estimate gives the solution to his 'problem of estimation', that is, likelihood/confidence regions must be formed using the observed likelihood function obtained from the (conditional) sampling distribution of the maximum likelihood estimate.

Pawitan (2001, p. 42) suggests the term *MLE-based regions (intervals)* to those formed from asymptotic normal distribution of the maximum likelihood estimate, that is, to *Wald confidence regions (intervals)* and uses the term *likelihood-based confidence regions (intervals)* for those obtained from likelihood function. In this article, this terminology is used.

The idea that by *Method of Maximum Likelihood* Fisher meant point estimation or confidence intervals based on an asymptotic normal distribution of MLE is in a complete contradiction with the fact that his *Theory of Statistical Estimation* was developed to justify the use of likelihood function in statistical inference. It is extremely odd to think that Fisher at the same time would have been promoting two solutions to his *Problem of Estimation*, especially if we take to account his life long emphasis of finite samples. Every time Fisher discussed the theory of statistical estimation when he started to discuss the case of finite samples he noted the unimportance of the large sample concepts. In Fisher (1956), for example, on page 147 he said that

> In fact, the asymptotic definition is satisfied by any statistic whatsoever applied to a finite sample, and is useless for the development of a theory of small samples.

and on page 159

> The theory of large samples can, however, never be more than a first step preliminary to the study of samples of finite size, . . . .

3.3. *Example.* Assume that the remission times in the control group of leukemia data may be modeled as a sample from some gamma distribution and that the standard deviation of the gamma distribution is the parameter of interest, that is,

$$\psi = \frac{\mu}{\sqrt{\lambda}},$$

where $\lambda$ and $\mu$ are the shape and mean parameters of the gamma distribution, respectively. Now the maximum likelihood estimate of $\psi$ is 6.763 and approximate 95%-level profile likelihood-based and delta-method confidence intervals for $\psi$ are $(4.634, 11.297)$ and $(3.829, 9.697)$, respectively. In a simulation of 10000 samples from gamma distribution with shape 1.642 and mean 8.667 the actual coverage probabilities were 0.938 and 0.837 for profile likelihood-based and delta-method intervals, respectively.

An explanation for the weak performance of delta-method in this case can be seen from Figure 5, which shows the log-likelihood function of $\psi$ for the actual data. The graph of log-likelihood function is asymmetric and so the profile likelihood-based confidence interval is also asymmetric with respect to the maximum likelihood estimate $\hat{\psi}$. The delta-method interval is forced to be symmetric and because of this it often misses the 'true' value of $\psi$ because of a too low upper limit. When applying delta-method, it is important to use an appropriate scale for the interest function by transforming first $\psi$ to $\lambda = h(\psi)$ so that the distribution of $\hat{\lambda} = h(\hat{\psi})$ is better approximated by normal distribution, especially with approximately constant variance. Then using the delta-method interval $(\hat{\lambda}^L, \hat{\lambda}^U)$ for $\lambda$ a better interval $(h^{-1}(\hat{\lambda}^L), h^{-1}(\hat{\lambda}^U))$ for $\psi$ is obtained. When applying the profile likelihood method, there is no need, however, to make any transformations, because of the invariance of profile likelihood-based intervals. In a sense, the profile likelihood method automatically chooses the best transformation and the user need not worry about that (Pawitan 2001, p. 47).

FIG 5. *Log-likelihood function of standard deviation of gamma distribution calculated from control group of leukamia data.*

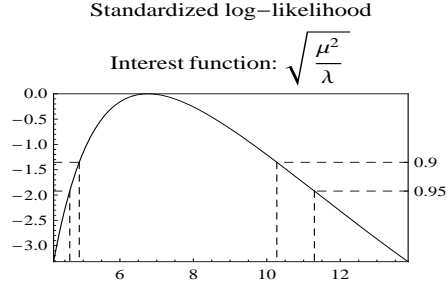

## 4. Current state of Fisher's theory.

4.1. *Fisher's problem of estimation after 1956.* Fisher lacked a systematic way of getting solution in case of general real-valued interest function $g(\omega)$ of multidimensional parameter. In Fisher (1922, p. 313) he briefly mentions the problem.

> There is one point, however, which may be briefly mentioned here in advance, as it has been the cause of some confusion. ...normal population of two correlated variates will usually require five parameters for its specification, the two means, the two standard deviations, and the correlation; of these often only the correlation is required, or if not alone of interest, it is discussed without reference to the other four quantities. In such cases an alteration has been made in what is, and what is not, relevent , and it is not surprising that certain small corrections should appear, or not, according as the other parameters of the hypothetical surface are not deemed relevant.

Fisher used marginal and conditional sampling distributions of carefully selected, but in anyway somewhat ad hoc, statistics and used likelihood functions obtained from these sampling distributions. These likelihood functions are special cases of so-called 'pseudo likelihood functions' (Pace and Salvan 1997, pp. 131-162). The term 'pseudo likelihood function' means a function of the interest parameter that is used instead of the original likelihood function. Marginal and conditional likelihood functions, are, however, genuine likelihood functions and the term 'pseudo' reflects the fact that they are used in place of original likelihood function and often do not contain all the information on the interest function included in statistical evidence. In addition there are other pseudo likelihood functions which are not likelihood functions at all. The prominent member of these other pseudo likelihood functions is *profile likelihood function*, which was explicitly introduced in Box and Cox (1962), but was implicitly used in connection of Wilk's (general) likelihood

ratio statistic

$$w(\psi; y_{obs}) = 2\{l(\hat{\omega}; y_{obs}) - l(\tilde{\omega}; y_{obs})\}$$

for the statistical hypothesis H : $g(\omega) = \psi$, where the function

$$l_P(\psi; y_{obs}) = l(\tilde{\omega}; y_{obs}) = \max_{\omega:g(\omega)=\psi} l(\omega; y_{obs})$$

is called *logarithmic profile likelihood function* or *profile log-likelihood function* of $\psi = g(\omega)$.

As discussed above the *MLE-based intervals* are against Fisher's ideology, because without appropriate transformations the intervals are not likelihood intervals. On the other hand *likelihood-based intervals* calculated from (profile) likelihood function are in accordance with Fisher's ideology.

4.2. *Profile likelihood function.* Let $y_{\text{obs}}$ be the observed response and $\mathcal{M} = \{p(y; \omega) : y \in \mathcal{Y}, \omega \in \Omega\}$ its statistical model. In most practical problems only part of the parameter vector or more generally the value of a given function of the parameter vector is of interest.

Let $g(\omega)$ be a given interest function with $q$-dimensional real vector $\psi$ as value. Then the function

$$(4.1) \qquad \begin{aligned} L_g(\psi; y_{\text{obs}}) &= \max_{\{\omega \in \Omega : g(\omega)=\psi\}} L(\omega; y_{\text{obs}}) \\ &= L(\tilde{\omega}_\psi; y_{\text{obs}}) \end{aligned}$$

is called the *profile likelihood function* of the interest function $g$ induced by the statistical evidence $(y_{\text{obs}}, \mathcal{M})$. The value $\tilde{\omega}_\psi$ of the parameter vector maximizes the likelihood function in the subset $\{\omega \in \Omega : g(\omega) = \psi\}$ of the parameter space. The function

$$(4.2) \qquad \begin{aligned} l_g(\psi; y_{\text{obs}}) &= \max_{\{\omega \in \Omega : g(\omega)=\psi\}} l(\omega; y_{\text{obs}}) \\ &= l(\tilde{\omega}_\psi; y_{\text{obs}}) \\ &= \ln(L_g(\psi; y_{\text{obs}})) \end{aligned}$$

is called the *logarithmic profile likelihood function* or *profile log-likelihood function* of the interest function $g$ induced by the statistical evidence $(y_{\text{obs}}, \mathcal{M})$. Furthermore functions

$$(4.3) \qquad R_g(\psi; y_{\text{obs}}) = \frac{L_g(\psi; y_{\text{obs}})}{L_g(\hat{\psi}; y_{\text{obs}})} = \frac{L_g(\psi; y_{\text{obs}})}{L(\hat{\omega}; y_{\text{obs}})}$$

and

$$(4.4) \qquad \begin{aligned} r_g(\psi; y_{\text{obs}}) &= l_g(\psi; y_{\text{obs}}) - l_g(\hat{\psi}; y_{\text{obs}}) \\ &= l_g(\psi; y_{\text{obs}}) - l(\hat{\omega}; y_{\text{obs}}) \end{aligned}$$

are called the *relative profile likelihood function* and the *logarithmic relative profile likelihood function* or *relative profile log-likelihood function* of interest function $g$, respectively. Because in this article no other likelihood functions except actual profile likelihood functions of various interest functions $g$ are considered the phrases (*relative*) *likelihood function* and (*relative*) *log-likelihood function* of $g$ are used.

When the interest function $g$ is real valued the parameter vectors $\tilde{\omega}_\psi$ form a curve in the parameter space. This curve is called the *profile curve* of the interest function $g$.

4.3. *Profile likelihood region and its uncertainty.* With help of the relative likelihood function of interest function $g$ one can construct the so called *profile likelihood regions*. The set

$$\{\psi : R_g(\psi; y_{\mathrm{obs}}) \geq c\} = \{\psi : r_g(\psi; y_{\mathrm{obs}}) \geq \ln(c)\}$$

of values of the interest function $g(\omega)$ is the $100c\%$ profile likelihood region. The value $\psi = g(\omega)$ of the parameter function does not belong to the $100c\%$ profile likelihood region if the response is such that

$$r_g(\psi; y) < \ln(c).$$

Probability of this event, calculated at a given value $\omega$ of the parameter vector, is used as a measure of uncertainty of the statement that the unknown value of the interest function belongs to the $100c\%$ profile likelihood region, provided that this probability has the same value or at least approximately the same value for all values of the parameter vector. One minus this probability is called the *confidence level* of the profile likelihood region and the region is called the (approximate) $(1 - \alpha)$-level *profile likelihood-based confidence region* for the interest function. Because under mild assumptions concerning the interest function $g(\omega)$ and statistical model the random variable $-2r_g(\psi; y)$ has approximately the $\chi^2[q]$-distribution, the set

$$\left\{\psi : r_g\left(\psi; y_{\mathrm{obs}}\right) \geq -\frac{\chi^2_{1-\alpha}[q]}{2}\right\} = \left\{\psi : -2r_g\left(\psi; y_{\mathrm{obs}}\right) \leq \chi^2_{1-\alpha}[q]\right\}$$

$$= \left\{\psi : l_g(\psi) \geq l_g(\hat{\psi}) - \frac{\chi^2_{1-\alpha}[q]}{2}\right\}$$

is the approximate $(1 - \alpha)$-level confidence region for the interest function $g(\omega)$ (Severini, 2000). In some cases the distribution of $-2r_g(\psi; y)$ is exactly the $\chi^2[q]$-distribution and then the set is exact confidence region. Sometimes the random variable $-2r_g(\psi; y)$ has some other known distribution, usually the $F$-distribution.

## 5. Case of real-valued interest function.

5.1. *Profile likelihood-based confidence interval.* For real valued interest functions the profile likelihood-based confidence regions are usually intervals and so those regions are called (approximate) $(1-\alpha)$-level *profile likelihood-based confidence intervals*. Thus the set

(5.1) $$\left\{ \psi : r_g \left( \psi ; y_{\text{obs}} \right) \geq -\frac{\chi^2_{1-\alpha}[1]}{2} \right\} = [\hat{\psi}_L, \hat{\psi}_U]$$

forms the (approximate) $(1-\alpha)$-level profile likelihood-based confidence interval for the interest function $\psi = g(\omega)$. The statistics $\hat{\psi}_L$ and $\hat{\psi}_U$ are called the *lower* and *upper limits* of the confidence interval, respectively.

5.2. *Calculation of profile likelihood-based confidence intervals.* If the (approximate) $(1-\alpha)$-level profile likelihood-based confidence region of the real valued interest function $\psi = g(\omega)$ is an interval, its end points $\hat{\psi}_L$ and $\hat{\psi}_U$ satisfy the relations $\hat{\psi}_L < \psi < \hat{\psi}_U$ and

(5.2) $$l_g(\hat{\psi}_L) = l_g(\hat{\psi}_U) = l_g(\hat{\psi}) - \frac{\chi^2_{1-\alpha}[1]}{2}.$$

Thus $\hat{\psi}_L$ and $\hat{\psi}_U$ are roots of the equation $l_g(\psi) = l_g(\hat{\psi}) - \frac{\chi^2_{1-\alpha}[1]}{2}$. Consequently most applications of the profile likelihood-based intervals have determined $\hat{\psi}_L$ and $\hat{\psi}_U$ using some iterative root finding method. This approach involves the solution of an optimization problem in every trial value and depending on the method also the calculation the derivatives of the logarithmic profile likelihood function at the same trial value.

The approximate $(1-\alpha)$-level profile likelihood confidence set for parameter function $\psi = g(\omega)$ satisfies the relation

$$\left\{ \psi : l_g(\psi) > l_g(\hat{\psi}) - \frac{\chi^2_{1-\alpha}[1]}{2} \right\} = \{g(\omega) : \omega \in \mathcal{R}_{1-\alpha}(y_{\text{obs}})\},$$

where

(5.3) $$\mathcal{R}_{1-\alpha}(y_{\text{obs}}) = \left\{ \omega : l(\omega) > l(\hat{\omega}) - \frac{\chi^2_{1-\alpha}[1]}{2} \right\}$$

is a likelihood region for the whole parameter vector. This result is true, because the real number $\psi$ belongs to the profile likelihood region, if and only if there exist a parameter vector $\omega^*$ such that $g(\omega^*) = \psi$ and

$$\sup_{\{\omega \in \Omega : g(\omega) = \psi\}} l(\omega) = l_g(\psi) > l(\omega^*) > l_g(\hat{\psi}) - \frac{\chi^2_{1-\alpha}[1]}{2}.$$

That is if and only if there exists a parameter vector $\omega^*$ belonging to the likelihood region $\mathcal{R}_{1-\alpha}(y_{\text{obs}})$ such that it satisfies the equation $g(\omega^*) = \psi$. So the number $\psi$ belongs to the profile likelihood region, if and only if it belongs to the set $\{g(\omega) : \omega \in \mathcal{R}_{1-\alpha}(y_{\text{obs}})\}$.

Assume now that the profile likelihood-based confidence region is an interval. Then the end points $\hat{\psi}_L$ and $\hat{\psi}_U$ satisfy the following relations

$$(5.4) \qquad \hat{\psi}_L = \inf_{\omega \in \mathcal{R}_{1-\alpha}(y_{\text{obs}})} g(\omega)$$

and

$$(5.5) \qquad \hat{\psi}_U = \sup_{\omega \in \mathcal{R}_{1-\alpha}(y_{\text{obs}})} g(\omega).$$

Let now the likelihood region be a bounded and connected subset of the parameter space. If the log-likelihood and the interest functions are continuous in the likelihood region, the latter with non-vanishing gradient, then the profile likelihood-based confidence region is an interval $(\hat{\psi}_L, \hat{\psi}_U)$ with

$$(5.6) \qquad \hat{\psi}_L = g(\tilde{\omega}_L) = \inf_{l(\omega) = l(\hat{\omega}) - \chi^2_\alpha[1]/2} g(\omega),$$

and

$$(5.7) \qquad \hat{\psi}_U = g(\tilde{\omega}_U) = \sup_{l(\omega) = l(\hat{\omega}) - \chi^2_\alpha[1]/2} g(\omega).$$

This follows from the assumptions, because they imply, that the set $\mathcal{R}_{1-\alpha}(y_{\text{obs}})$ is open, connected, and bounded subset of the $d$-dimensional space. Thus the closure of $\mathcal{R}_{1-\alpha}(y_{\text{obs}})$ is a closed, connected, and bounded set. Form the assumptions concerning $g$ it follows that it attains its infimum and supremum on the boundary of $\mathcal{R}_{1-\alpha}(y_{\text{obs}})$ and takes every value between infimum and supremum somewhere in $\mathcal{R}_{1-\alpha}(y_{\text{obs}})$.

Under the above assumptions the solutions of the following constrained minimization (maximization) problem

$$(5.8) \qquad \min(\max) g(\omega)$$

with constraint

$$(5.9) \qquad l(\omega; y_{\text{obs}}) = l(\hat{\omega}; y_{\text{obs}}) - \frac{\chi^2_\alpha[1]}{2}.$$

gives the lower (upper) limit point of the profile likelihood-based confidence interval. This problem is rarely explicitly solvable and requires use of some kind of iteration (Uusipaikka 1996, Virtanen and Uusipaikka 2008).

5.3. *Modifications to obtain better coverage probabilities.* There exists various ways to modify likelihood quantities so that the confidence intervals obtained from these modifications have better coverage properties. Bartlett-correction (Bartlett 1937, 1938; Box 1946) is used to modify the likelihood ratio statistic, Lugannani-Rice (1980) modifies the signed root deviance, that is, the signed square root of the likelihood ratio statistic, and many people have suggested various modifications of the profile likelihood function. References to articles discussing these can be found, for example, in Barndorff-Nielsen and Cox (1994) and Severini (2000).

All these modifications have been used in the way that even in the case of a one-dimensional parameter the obtained intervals are not likelihood intervals. Because likelihood is the primary concept and coverage probability a secondary one this is unfortunate. It should be possible to use these modifications so that the obtained intervals would still be likelihood intervals, but with better repeated sampling properties.

5.4. *Example.* Assume that the times of remission can be considered as samples from two exponential distributions. Statistical model under these assumptions is

$$\mathcal{M} = \text{IndependenceModel}[\{\mathcal{M}_1, \mathcal{M}_2\}],$$

where

$$\mathcal{M}_1 = \text{SamplingModel}[\text{ExponentialModel}[\mu_1], 21],$$

and

$$\mathcal{M}_2 = \text{SamplingModel}[\text{ExponentialModel}[\mu_2], 21]$$

are models for remission times of *Drug* and *Control* groups, respectively.

Of interest in addition to model parameters might be the difference of means or their ratio, that is, parameter functions

$$\psi_1 = \mu_1 - \mu_2$$

and

$$\psi_2 = \mu_1/\mu_2,$$

but better characteristic for describing the difference of distributions might be the probability that random observation from the first distribution would be greater than random statistically independent observation from second distribution. This interest function is related to the so-called *Receiver Operating Characteristic* (ROC) curve, namely, it is area under this curve (AUC). Under current assumptions, this interest function has the form

$$\psi_3 = \frac{\mu_1}{\mu_1 + \mu_2}.$$

TABLE 3
*Estimates and approximate 0.95-level confidence intervals of interest functions.*

| Parameter function | Estimate | Lower limit | Upper limit |
|---|---|---|---|
| $\mu_1$ | 39.889 | 22.127 | 83.013 |
| $\mu_2$ | 8.667 | 5.814 | 13.735 |
| $\psi_1$ | 31.222 | 12.784 | 74.444 |
| $\psi_2$ | 4.603 | 2.174 | 10.583 |
| $\psi_3$ | 0.822 | 0.685 | 0.914 |
| $\psi_4$ | -2.820 | -8.677 | -0.634 |

FIG 6. *Plot of log-likelihood of* $\psi_2$

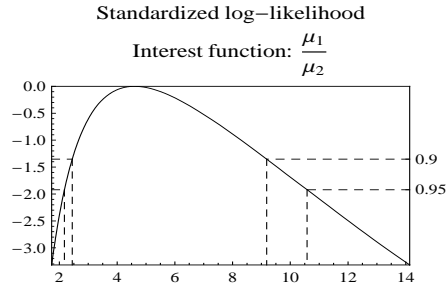

Another parameter function describing the difference of the distributions is the so-called *Kullback-Leibler divergence*, which under current assumptions is

$$\psi_4 = 2 - \frac{\mu_2}{\mu_1} - \frac{\mu_1}{\mu_2}.$$

Table 3 gives maximum likelihood estimates and approximate 0.95-level profile likelihood-based confidence intervals for parameters and interest functions. Figures 6 and 7 show graphs of log-likelihoods of $\psi_2$ and $\psi_3$, respectively.

**6. Conclusions.** In this article a possible, in author's mind highly plausible and well justified, interpretation to Fisher's meaning for the word 'estimate' has been given. According to this interpretation estimate is a statistic taking values in the parameter space such that the observed likelihood function obtained from the sampling distribution of this estimate is used to construct likelihood intervals for the parameter. Fisher's *Theory of Statistical Estimation* is a chain of ideas which Fisher developed to justify his solution to the problem of estimation. This solution is Method of Maximum Likelihood in which likelihood intervals are produced from the observed likelihood function obtained from the (conditional) sampling distribution of the maximum likelihood estimate. As shown above the end result of Fisher's justification is that likelihood intervals must be based on the original likelihood function obtained from statistical evidence so that all relevant information will be used.

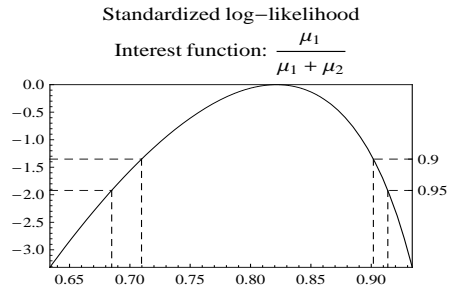

FIG 7. *Plot of log-likelihood of $\psi_3$.*

MLE-based confidence intervals calculated from the asymptotic normal distribution of the maximum likelihood estimate using delta-method is the current dominant approach. This approach was shown to involve Fisher's ideas but to be against his ideology. They are against Fisher's ideology because generally MLE-based confidence intervals are not likelihood intervals. In addition these intervals are not invariant, may contain impossible values, and have poor coverage probabilities.

Fisher did not have a systematic solution of the problem of estimation when a real-valued parameter function of a higher dimensional parameter is of interest. During last fifty years, various solutions based on pseudo likelihood functions have been suggested. The one that has generated most research and applications is likelihood-based inference using profile likelihood function. It was shown that this gives general solution to Fisher's problem of estimation.

A new method of calculation of profile likelihood-based intervals was considered. This method is a simple powerful method that can be used for general smooth interest functions in general statistical models. Therefore there is no theoretical and practical reason to calculate MLE-based confidence intervals instead of profile likelihood-based intervals.

In conclusion, even though it seems to be very common opinion, for Fisher the method of maximum likelihood did not mean the usage of MLE as a point estimate or usage of the asymptotic normal distribution of MLE to construct confidence intervals. Instead it meant the construction of likelihood intervals or regions from original likelihood function or some pseudo likelihood function and so interpreted Fisher's method of maximum likelihood is very much like the method of support discussed by A. W. F. Edwards in his book *Likelihood* (1972).

## References.

DEPARTMENT OF STATISTICS
UNIVERSITY OF TURKU
FIN-20014 TURKU
FINLAND E-MAIL: esa.uusipaikka@utu.fi